\documentclass[]{imsart}
\usepackage{amsthm,amsmath}
\usepackage[colorlinks,citecolor=blue,urlcolor=blue]{hyperref}

\usepackage{amsmath,amsfonts,amssymb,amsthm}

\usepackage{enumitem}

\RequirePackage{amsthm,amsmath,amsfonts,amssymb}
\RequirePackage[numbers]{natbib}
\RequirePackage[colorlinks,citecolor=blue,urlcolor=blue]{hyperref}
\RequirePackage{graphicx}
\usepackage{lineno,hyperref}

\startlocaldefs
\modulolinenumbers[5]

\usepackage[utf8]{inputenc}
\usepackage[english]{babel}

\usepackage[margin=1in]{geometry}
\usepackage{bm}
\usepackage{graphicx}
\usepackage{amssymb,amsmath,amsthm,mathrsfs}
\usepackage{epstopdf}
\usepackage{algorithm}
\usepackage{amsfonts,delarray}
\usepackage{algorithm,algcompatible,amsmath}
\usepackage{rotating,color}
\usepackage{subfigure}
\usepackage{multirow}
\usepackage{titlesec,textcase}
\usepackage{longtable}
\usepackage{bbm}
\usepackage{rotating}

\newtheorem{Theorem}{Theorem}[section]
\newtheorem{Lemma}[Theorem]{Lemma}

\newtheorem{Proposition}[Theorem]{Proposition}

\theoremstyle{definition}

\RequirePackage[normalem]{ulem} 

\definecolor{rp}{RGB}{83,54,106}

\def\boxit#1{\vbox{\hrule\hbox{\vrule\kern6pt\vbox{\kern6pt#1\kern6pt}\kern6pt\vrule}\hrule}}

\endlocaldefs

\begin{document}
\begin{frontmatter}
\title{Statistical Limits for Testing Correlation of Hypergraphs}

\runtitle{Hypergraph Correlation Testing}
\runauthor{Mingao Yuan and Zuofeng Shang}
\begin{aug}
\author[A]{\fnms{Mingao} \snm{Yuan}\ead[label=e1]{mingao.yuan@ndsu.edu}}
\and
\author[B]{\fnms{Zuofeng} \snm{Shang}\ead[label=e2]{zshang@njit.edu}}
\address[A]{Department of Statistics, 
North Dakota State University;
\printead{e1}}

\address[B]{Department of Mathematical Sciences,
New Jersey Institute of Technology;
\printead{e2}}
\end{aug}

\begin{abstract}
In this paper, we consider the hypothesis testing of correlation between two $m$-uniform hypergraphs on $n$ unlabelled nodes. Under the null hypothesis, the hypergraphs are independent, while under the alternative hypothesis, the hyperdges have the same marginal distributions as in the null hypothesis but are correlated after some unknown node permutation. We focus on two scenarios: the hypergraphs are generated from the Gaussian-Wigner model and the dense Erd\"{o}s-R\'{e}nyi model. We derive the sharp information-theoretic testing threshold. Above the threshold, there exists a powerful test to distinguish the alternative hypothesis from the null hypothesis. Below the threshold, the alternative hypothesis and the null hypothesis are not distinguishable. The threshold involves $m$ and decreases as $m$ gets larger. This indicates testing correlation of hypergraphs ($m\geq3$) becomes easier than testing correlation of graphs ($m=2$).
\end{abstract}



\end{frontmatter}

\section{Introduction}
\label{S:1}

Graph matching is a fundamental problem in network data analysis. It refers to the problem of identifying a mapping between the nodes of two graphs that preserves as much as possible the relationships between nodes. Graph matching is a powerful technique and  is widely used in a variety of scientific fields. For instance, in shape matching and object recognition, graph matching is used to find the correspondence between object graph and its feature graph(\cite{BBM05,CL12}); in social network analysis, graph matching identifies all the accounts belonging to the same individual (\cite{KL14}); in computational biology,  graph matching can be applied to match brain-graphs(\cite{VCLPKH15}). Graph matching problem is NP hard in the worst case and various algorithms have been developed to recover the latent mapping (\cite{CSS06,VCLPKH15,KL14,CL12,BCLSS19,BBM05,DMWX21, YXL21}). In practice, whether there exists a meaningful matching between two graphs is unknown. To solve this issue, \cite{BCLSS19,WXY21,CWXY21} initiate the study of testing the correlation of two graphs. Especially, \cite{WXY21} derives the sharp information-theoretic threshold for testing correlated Gaussian-Wigner graphs and dense Erd\"{o}s-R\'{e}nyi graphs and \cite{CWXY21} propose a test procedure with polynomial-time complexity.

Many complex networks in the real world can be formulated as hypergraphs. Unlike ordinary graphs where the data structure is typically unique, e.g., edges only contain two vertices, hypergraphs demonstrate a number of possibly overlapping data structures so that an edge may contain arbitrarily many vertices. For instance, in coauthorship networks
(\cite{ER05,OGM17, RDP04,N01}), an edge represents a group of arbitrarily many coauthors; in folksonomy network, an edge may represent a triple (user, resource, annotation) structure  (\cite{GZCN09}); in login network an edge may represent a (user, remote host, login time, logout time) structure (\cite{GD14}). Recently, there is increasing interest in hypergraph matching problem, that is,  to establish the correspondence between nodes of two unlabelled hypergraphs (\cite{ZS08,DBKP11,LCL11,NTGH16,PPH13}). In this paper, we study the hypothesis testing of correlation for hypergraphs and characterize how the sharp testing threshold in \cite{WXY21} varies in hypergraph. 

An \textit{undirected} $m$-uniform hypergraph is a pair $\mathcal{H}_m=([n],\mathcal{E})$ in which $[N]:=\{1,2,\dots,n\}$ 
is a vertex set and $\mathcal{E}$ is a set of hyperedges. Each hyperedge in $\mathcal{E}$ consists of exactly $m$ vertices in $[n]$. 
The corresponding adjacency tensor is an $m$-dimensional symmetric array $A\in(B^n)^{\otimes m}$ 
satisfying $A_{i_1i_2\ldots i_m}\in B$ for $1\leq i_1<i_2<\dots<i_m\leq n$, in which $B\subset\mathbb{R}$. 
Here, symmetry means that $A_{i_1i_2\ldots i_m}=A_{j_1j_2\ldots j_m}$ whenever $i_1,i_2,\ldots,i_m$ is a permutation of $j_1,j_2,\ldots,j_m$. 
If $|\{i_1,i_2,\ldots,i_m\}|\leq m-1$, then $A_{i_1i_2\ldots i_m}=0$,
i.e., no self-loops are allowed. 
In particular, $B=\{0,1\}$ corresponds to binary hypergraphs.
The general $B$ corresponds to weighted hypergraphs. For convenience, we also denote the hypergraph $\mathcal{H}_m=([n],\mathcal{E})$ as $\mathcal{H}_m=([n],A)$.

Let $P_n$ be the permutation group on $[n]$. 
Two hypergraphs $\mathcal{H}_{m,1}=([n],A_1)$ and $\mathcal{H}_{m,2}=([n],A_2)$ are said to be isomorphic, denoted as $\mathcal{H}_{m,1}\cong\mathcal{H}_{m,2}$ if there is a permutation $\pi\in P_n$ such that $A_{1,i_1i_2\dots i_m}=A_{2,\pi_{i_1}\pi_{i_2}\dots \pi_{i_m}}$ for all $1\leq i_1<i_2<\dots<i_m\leq n$. Clearly isomorphism defines an equivalence relation and denote the equivalence class of $\mathcal{H}_{m,1}$ as $\overline{\mathcal{H}}_{m,1}$. Each hypergraph $\mathcal{H}_{m}\in \overline{\mathcal{H}}_{m,1}$ is called an unlabelled hypergraph of $\mathcal{H}_{m,1}$.

For two hypergraphs $\mathcal{H}_{m,1}=([n],A_1)$ and $\mathcal{H}_{m,2}=([n],A_2)$, suppose $(A_{1,i_1i_2\dots i_m},A_{2,i_1i_2\dots i_m})$, $(1\leq i_1<i_2<\dots<i_m\leq n)$ are independently and identically distributed random variables with $A_{1,i_1i_2\dots i_m}$ and $A_{2,i_1i_2\dots i_m}$ sharing the same marginal distribution. Given two unlabelled hypergraphs (random sample) $\widetilde{A}_1\in\overline{\mathcal{H}}_{m,1}$ and $\widetilde{A}_2\in\overline{\mathcal{H}}_{m,2}$, our purpose is to test the following hypergraph correlation hypothesis.
\begin{eqnarray}\nonumber
    &&H_0:\ A_{1,i_1i_2\dots i_m}\ and\ A_{2,i_1i_2\dots i_m}\ are\ independent, \ \ 1\leq i_1<i_2<\dots<i_m\leq n; \\ \label{hypothesis}
  &&  
     H_1: \ A_{1,i_1i_2\dots i_m}\ and\ A_{2,i_1i_2\dots i_m}\ are\ correlated, \ \ \ \ \ 1\leq i_1<i_2<\dots<i_m\leq n.
\end{eqnarray}
When $m=2$, (\ref{hypothesis}) is just the graph correlation hypothesis testing problem studied in \cite{BCLSS19,WXY21,CWXY21}. It is not immediately clear what role $m\geq3$ plays in the hypothesis testing problem (\ref{hypothesis}). This motivates us to study (\ref{hypothesis}) for general $m\geq2$.

In this paper, we focus on two scenarios.
\begin{itemize}
\item[] (I) {\bf Gaussian-Wigner hypergraph}: For all $1\leq i_1<i_2<\dots<i_m\leq n$, $A_{1,i_1i_2\dots i_m}$ and $A_{2,i_1i_2\dots i_m}$ follow the bivariate normal distribution with mean zero, variance one and correlation coefficient $\rho\in[0,1]$. Then (\ref{hypothesis}) is simplified to $H_0:\rho=0,\ v.s.\ H_1: \rho>0$. 

\item[] (II) {\bf Erd\"{o}s-R\'{e}nyi hypergraph}: Let $\mathcal{H}_{m}$ and $\mathcal{H}_{m}^{\prime}$ be independent random Erd\"{o}s-R\'{e}nyi $m$-uniform hypergraphs with hyperedge probability $p\in[0,1]$. The we can restate (\ref{hypothesis}) as follows: $H_0$ is equivalent to that $\mathcal{H}_{m,1}$ and $\mathcal{H}_{m,2}$ are generated from $\mathcal{H}_{m}$ and $\mathcal{H}_{m}^{\prime}$ respectively by keeping each hyperedge independently with probability $s\in [0,1]$; $H_1$ is equivalent to that $\mathcal{H}_{m,1}$ and $\mathcal{H}_{m,2}$ are similarly generated from the same hypergraph $\mathcal{H}_{m}$. In this case, the correlation between $(A_{1,i_1i_2\dots i_m},A_{2,i_1i_2\dots i_m})$, $(1\leq i_1<i_2<\dots<i_m\leq n)$ under $H_1$ is $\rho=\frac{s(1-p)}{1-ps}$.
\end{itemize}
We shall use the total variation distance to measure the difference between $H_1$ and $H_0$.
The total variation distance between two probability measures $P,Q$ on a sigma-algebra $\mathcal{F}$ of subsets of the sample space $\Omega$ is defined as
\[TV(P,Q)=\sup_{E\in\mathcal{F}}|P(E)-Q(E)|.\]
Let $P,Q$ be probability measures under $H_0, H_1$ respectively.
Then $H_0$ and $H_1$ are said to be indistinguishable if $TV(P,Q)=o(1)$ and distinguishable if $TV(P,Q)=1+o(1)$.

In this paper, we adopt the Bachmann-Landau notation $o(1),O(1)$. For two positive sequences $a_n,b_n$, denote $a_n\asymp b_n$ or $a_n=\Theta(b_n)$ if $ 0<c_1\leq\frac{a_n}{b_n}\leq c_2<\infty$ for constants $c_1,c_2$. Denote $a_n\gg b_n$ or $b_n\ll a_n$ if $\lim_{n\rightarrow\infty}\frac{a_n}{b_n}=\infty$. We write $a_n=\Omega(b_n)$ if $a_n\geq cb_n$ for a constant $c>0$. $I[E]$ denotes the indicator function of event $E$.

The rest of the paper is organized as follows. In section 2, we present the main result and related proof for Gaussian-Wigner Model. Section 3 provides the main result and proof for Erd\"{o}s-R\'{e}nyi Model. Some necessary lemmas are given in section 4.

\section{Gaussian-Wigner Hypergraph}
In this section, we study the hypergraph correlation test problem under the Gaussian-Wigner model. Denote $\pi\sim Unif(P_n)$ if $\pi$ is uniformly and randomly selected from $P_n$. In this case, the hypothesis (\ref{hypothesis}) is reformulated as follows. 
\begin{eqnarray}\nonumber
H_0:
\begin{pmatrix}A_{1,i_1i_2\dots i_m}\\
A_{2,i_1i_2\dots i_m}
\end{pmatrix} & \overset{i.i.d.}{\sim} & N\left[\left(\begin{array}{c}
0\\
0
\end{array}\right),\left(\begin{array}{ccc}
 1 & 0\\
 0 & 1\\
\end{array}\right)\right],\\ \label{gaussianhyp}
H_1:
\begin{pmatrix}A_{1,i_1i_2\dots i_m}\\
A_{2,\pi_{i_1}\pi_{i_2}\dots \pi_{i_m}}
\end{pmatrix} & \overset{i.i.d.}{\sim} & N\left[\left(\begin{array}{c}
0\\
0
\end{array}\right),\left(\begin{array}{ccc}
 1 & \rho\\
 \rho & 1\\
\end{array}\right)\right],\ conditional\ on\ \pi\sim Unif(P_n).
\end{eqnarray}
When $m=2$, the Gaussian-Wigner model is proposed in \cite{DMWX21} and studied in \cite{GLM19,FMWX19,WXY21}. The following result provides the sharp information-theoretic threshold for hypothesis testing problem (\ref{gaussianhyp}).
\begin{Theorem}[Gaussian-Wigner hypergraph]\label{thm:1}
Let $m\geq2$ be any fixed integer. Then $H_0$ and $H_1$ in (\ref{gaussianhyp}) are distinguishable if
\[\rho^2\geq\frac{2n\log n}{\binom{n}{m}}.\]
$H_0$ and $H_1$ in (\ref{gaussianhyp}) are indistinguishable if
\begin{equation}\label{rhos}
\rho^2<\frac{(1-\epsilon)2n\log n}{\binom{n}{m}},
\end{equation}
for any constant $\epsilon>0$.
\end{Theorem}
 
 For Gaussian-Wigner model, a phase transition phenomenon occurs at the threshold $\frac{2n\log n}{\binom{n}{m}}$: $H_1$ and $H_0$ are distinguishable if and only if the correlation is above the threshold. Note that the threshold decreases at rate $\frac{\log n}{n^{m-1}}$ as a function of $m$. This indicates that testing correlated Gaussian-Wigner hypergraphs ($m\geq3$) is easier than testing correlated Gaussian-Wigner graphs (see result for $m=2$ in \cite{WXY21}).

\begin{proof}[Proof of Theorem \ref{thm:1}] 
{\bf (Positive result).}  We shall construct a powerful test statistic based on the maximum likelihood method. Since the testing problem is easier for larger $\rho^2$, then we can assume $\rho^2=\frac{2n\log n}{\binom{n}{m}}$. For convenience, let $t_n=\rho \binom{n}{m}-\sqrt{\binom{n}{m}}n^{0.25}$.

Let $\pi$ be a uniformly and randomly selected permutation on $[n]$ such that $A_{1,i_1i_2\dots i_m}$ and $A_{2,\pi_{i_1}\pi_{i_2}\dots \pi_{i_m}}$ follow the bivariate normal distribution with mean zero, variance one and correlation coefficient $\rho\in[0,1]$. 

Under $H_1$, the likelihood ratio given $\pi$ is equal to
\begin{eqnarray}\nonumber
\frac{Q(A_1,A_2|\pi)}{P(A_1,A_2)}&=&\frac{1}{\sqrt{1-\rho^2}^{\binom{n}{m}}}\exp\left\{-\frac{\rho^2}{2(1-\rho^2)}\sum_{1\leq i_1<\dots<i_m\leq n}(A_{1,i_1i_2\dots i_m}^2+A_{2,\pi_{i_1}\pi_{i_2}\dots \pi_{i_m}}^2)\right\}\\ \label{likelihoodg}
&&\times\exp\left\{\frac{\rho}{1-\rho^2}\sum_{1\leq i_1<\dots<i_m\leq n}A_{1,i_1i_2\dots i_m}A_{2,\pi_{i_1}\pi_{i_2}\dots \pi_{i_m}}\right\}.
\end{eqnarray}
Hence, to maximize the likelihood ratio with respect to $\pi$ is equivalent to maximizing $T(\pi)$ given by
\[
T(\pi)=\sum_{1\leq i_1<\dots<i_m\leq n}A_{1,i_1i_2\dots i_m}A_{2,\pi_{i_1}\pi_{i_2}\dots \pi_{i_m}}.
\]
Then we define the test statistic as $T_n=\max_{\pi}T(\pi)$.

Under the alternative hypothesis, we shall show $\mathbb{P}(T_n\geq t_n)=1+o(1)$. By the Hanson-Wright inequality in Lemma \ref{hansonwright}, it is easy to verify that
\begin{eqnarray*}
\mathbb{P}(T_n\leq t_n)\leq \mathbb{P}(T(\pi)\leq t_n)\leq e^{-cn^{\frac{m+0.5}{2}}}+e^{-c\sqrt{n}},
\end{eqnarray*}
for some constant $c>0$. Then $\mathbb{P}(T_n\leq t_n)=o(1)$ and hence $\mathbb{P}(T_n\geq t_n)=1+o(1)$.

Under the null hypothesis, we show $\mathbb{P}(T_n\geq t_n)=o(1)$. Note that $A_{1,i_1i_2\dots i_m}$ and $A_{2,\pi_{i_1}\pi_{i_2}\dots \pi_{i_m}}$ are independent  for any $\pi$ and they follow the standard normal distribution. For  $\lambda:=\frac{t_n}{\binom{n}{m}}=o(1)$, the Chernoff bound in Lemma \ref{chernoff} yields
\begin{eqnarray*}
\mathbb{P}(T(\pi)\geq t_n)&=&\mathbb{P}(e^{T(\pi)}\geq e^{t_n})\leq \exp\left\{-\lambda t_n-\frac{\binom{n}{m}}{2}\log(1-\lambda^2)\right\}\\
&=&\exp\left\{-2n\log n-n^{0.5}+2\sqrt{2n\log n}n^{0.25}+n\log n+\frac{n^{0.5}}{2}+o(n)\right\}.
\end{eqnarray*}
Note that $n!\leq en^{n+0.5}e^{-n}$. Then by the union bound, it follows that
\[\mathbb{P}(T_n\geq t_n)\leq n!\mathbb{P}(T(\pi)\geq t_n)=\exp(-n+o(n))=o(1).\]
Then the proof is complete.
\end{proof}

\begin{proof}[Proof of Theorem \ref{thm:1}] ({\bf Negative result}).  To prove the negative result, it suffices to prove the second moment of the likelihood ratio under $H_0$ converges to one, that is,
\[\mathbb{E}\left[\left(\frac{Q(A_1,A_2)}{P(A_1,A_2)}\right)^2\right]\leq 1+o(1),\]
under $H_0$. The details are given in the following Proposition \ref{prho} and Proposition \ref{prho2}.
\end{proof}

Before presenting Proposition \ref{prho} and Proposition \ref{prho2}, we provide some basic facts about permutation.
Each permutation $\pi\in P_n$ can be decomposed into product of disjoint cycles. Each cycle forms an orbit of any element in the cycle. Let $K_m$ be the complete $m$-uniform hypergraph on $[n]$. Then $\pi$ induces a permutation $\pi^K$ on the hyperedge set of $K_m$ by 
\[\pi^K(i_1,i_2,\dots,i_m)=(\pi_{i_1}, \pi_{i_2},\dots,\pi_{i_n}),\ \ i_1<i_2<\dots<i_m.\]
We call $\pi$ node permutation and $\pi^K$ hyperedge permutation. Let $n_k$ denote the number of cycles (orbits) in $\pi$ with length $k$ and $N_k$ the number of hyperedge cycles (hyperedge orbits) with length $k$. Note that $N_k$ can be expressed as a function of $n_t,(t\leq k)$. For example, let $m=3$. Then $N_1=\binom{n_1}{3}+n_1n_2+n_3$.

\begin{Proposition}\label{prho}
For any fixed integer $m\geq2$, if $\rho^2<\frac{(1-\epsilon)n\log n}{\binom{n}{m}}$ for any  constant $\epsilon\in[0,1)$, then $H_0$ and $H_1$ are indistinguishable for both Gaussian Wigner model and Erdos-Renyi model.
\end{Proposition}

\begin{proof}[Proof of Proposition \ref{prho}] 
We only need to focus on $m\geq3$, since the result for $m=2$ is given in \cite{WXY21}.
Denote $\tilde{\pi}$ be an independent copy of $\pi$. Firstly, we consider Gaussian Wigner model. Define
\[L_1(A_{1,i_1i_2\dots i_m},A_{2,\pi_{i_1}\pi_{i_2}\dots \pi_{i_m}})=\frac{1}{\sqrt{1-\rho^2}}\exp\left\{\frac{-\rho^2(A_{1,i_1i_2\dots i_m}^2+A_{2,\pi_{i_1}\pi_{i_2}\dots \pi_{i_m}}^2)+2\rho A_{1,i_1i_2\dots i_m}A_{2,\pi_{i_1}\pi_{i_2}\dots \pi_{i_m}}}{2(1-\rho^2)}\right\},
\]
and
\[L_{i_1i_2\dots i_m}=L_1(A_{1,i_1i_2\dots i_m},A_{2,\pi_{i_1}\pi_{i_2}\dots \pi_{i_m}})L_1(A_{1,i_1i_2\dots i_m},A_{2,\tilde{\pi}_{i_1}\tilde{\pi}_{i_2}\dots \tilde{\pi}_{i_m}}).\]

By (\ref{likelihoodg}), the second moment of the likelihood ratio under $H_0$ is equal to
\begin{eqnarray}\label{e2}
\mathbb{E}\left[\left(\frac{Q(A_1,A_2)}{P(A_1,A_2)}\right)^2\right]=\mathbb{E}_{\pi,\tilde{\pi}}\left(\mathbb{E}\left[\frac{Q(A_1,A_2|\pi)}{P(A_1,A_2)}\frac{Q(A_1,A_2|\tilde{\pi})}{P(A_1,A_2)}\right]\right)&=&\mathbb{E}_{\pi,\tilde{\pi}}\left(\mathbb{E}\prod_{1\leq i_1<\dots<i_m\leq n}L_{i_1i_2\dots i_m}\right).
\end{eqnarray}
Denote $\sigma=\pi^{-1}\circ \tilde{\pi}$. For a hyperedge orbit $O$ induced by $\sigma$, define
\[L_O=\prod_{\{i_1,\dots,i_m\}\in O}L_{i_1i_2\dots i_m}.\]
Since $\tilde{\pi}(e)=\pi\circ \sigma (e)$ for any hyperedge $e$, then $L_O$ only depends on $A_{1,e}, A_{2,\pi_e}$ for $e\in O$.

Let $\mathcal{O}$ be the set of hyperedge orbits of $\sigma$. Note that the hyperedge orbits are mutually disjoint and $A_{1,i_1i_2\dots i_m}$ and $A_{2,i_1i_2\dots i_m}$ are i.i.d. under $H_0$. Then by (\ref{e2}), we have
\begin{equation}\label{es2}
\mathbb{E}\left[\left(\frac{Q(A_1,A_2)}{P(A_1,A_2)}\right)^2\right]=\mathbb{E}_{\pi,\tilde{\pi}}\left(\prod_{O\in\mathcal{O}}\mathbb{E}(L_{O})\right)=\mathbb{E}_{\pi,\tilde{\pi}}\left[\prod_{k=1}^{\binom{n}{m}}\left(\frac{1}{1-\rho^{2k}}\right)^{N_k}\right],
\end{equation}
where the second equality follows from Proposition 1 in \cite{WXY21} and $N_k$ is the number of hyperedge orbits with length $k$.

Note that $\sum_{k=2}^{\binom{n}{m}}N_k\leq n^m$. According to (\ref{rhos}),  $\rho^4n^m=O\left(\frac{\log n}{n^{m-2}}\right)=o(1)$
for $m\geq3$. Consequently,
\[\prod_{k=2}^{\binom{n}{m}}\left(\frac{1}{1-\rho^{2k}}\right)^{N_k}\leq \left(\frac{1}{1-\rho^{4}}\right)^{\binom{n}{m}}\leq \exp\left(\frac{n^m\rho^4}{1-\rho^4}\right)=1+o(1).\]
Then
\begin{equation}\label{es3}
\mathbb{E}\left[\left(\frac{Q(A_1,A_2)}{P(A_1,A_2)}\right)^2\right]\leq(1+o(1))\mathbb{E}_{\pi,\tilde{\pi}}\left[\left(\frac{1}{1-\rho^{2}}\right)^{N_1}\right]\leq (1+o(1))\mathbb{E}_{\pi,\tilde{\pi}}\left[\exp\left(\frac{N_1\rho^2}{1-\rho^{2}}\right)\right]\leq1+o(1).
\end{equation}
where the last step follows from
the following Lemma \ref{enrho}. Then the proof is complete for Gaussian Wigner model.

For Erdos-Renyi model, by a similar argument and using Proposition 1 in \cite{WXY21}, we have 
\begin{equation*}\label{es2}
\mathbb{E}\left[\left(\frac{Q(A_1,A_2)}{P(A_1,A_2)}\right)^2\right]=\mathbb{E}_{\pi,\tilde{\pi}}\left[\prod_{k=1}^{\binom{n}{m}}\left(1+\rho^{2k}\right)^{N_k}\right].
\end{equation*}
By the condition $\rho^2<\frac{(1-\epsilon)n\log n}{\binom{n}{m}}$ and $m\geq3$, it follows that
\[\prod_{k=2}^{\binom{n}{m}}\left(1+\rho^{2k}\right)^{N_k}\leq \left(1+\rho^{4}\right)^{\binom{n}{m}}\leq \exp\left(n^m\rho^4\right)=1+o(1).\]
Hence
\begin{equation*}\label{es2}
\mathbb{E}\left[\left(\frac{Q(A_1,A_2)}{P(A_1,A_2)}\right)^2\right]\leq (1+o(1))\mathbb{E}_{\pi,\tilde{\pi}}\left[\left(1+\rho^{2}\right)^{N_1}\right]\leq (1+o(1))\mathbb{E}_{\pi,\tilde{\pi}}\left[\exp\left(N_1\rho^{2}\right)\right]\leq 1+o(1).
\end{equation*}
Here the last inequality follows from Lemma \ref{enrho}.
\end{proof}

\begin{Lemma}\label{enrho}
Let $N_1$ be the number of hyperedge orbits of $\sigma=\pi^{-1}\circ\tilde{\pi}$ with length one. If $\rho^2<\frac{(1-\epsilon)n\log n}{\binom{n}{m}}$ for any positive constant $\epsilon$, then
\[\mathbb{E}_{\pi,\tilde{\pi}}\left[\exp\left(\frac{N_1\rho^2}{1-\rho^{2}}\right)\right]\leq1+o(1).\]
\end{Lemma}

\begin{proof}[Proof of Lemma \ref{enrho}] 
Let $n_k$ be the number of $k$-nodes cycles of permutation $\sigma$. Since the cycles of $\sigma$ are disjoint, then $n_k\leq n$. Note that 1-hyperedge orbit is just a single hyperedge and this hyperedge can only involve nodes in $k$-nodes cycles with $k\leq m$. Hence,
$N_1=R(n_1,n_2,\dots,n_m)$,
where $R(n_1,n_2,\dots,n_m)$ is a polynomial in $n_1,n_2,\dots,n_m$. 
If a hyperedge contains a $k$-node cycle, then we only need to select $m-k$ nodes to form a hyperedge. Hence,
any terms in $R(n_1,n_2,\dots,n_m)$ involving $k$-node cycles are bounded by $n_kn^{m-k}=O(n^{m-k+1})$. Since $\rho^2<\frac{(1-\epsilon)n\log n}{\binom{n}{m}}$, then $\rho^2n^{m-k+1}=O\left(\frac{\log n}{n^{k-2}}\right)=o(1)$ for $k\geq3$.
If a term in $R(n_1,n_2,\dots,n_m)$ contains $n_2^k$, then it is bounded by $\rho^2n_2^kn^{m-2k}=O\left(\frac{\log n}{n^{k-1}}\right)=o(1)$ for $k\geq2$. Hence, we have
\[\rho^2N_1=\rho^2\left[\binom{n_1}{m}+n_2\binom{n_1}{m-2}\right]+o(1).\]
Then
\begin{eqnarray*}
\mathbb{E}_{\pi,\tilde{\pi}}\left[\exp\left(\frac{N_1\rho^2}{1-\rho^{2}}\right)\right]&=&(1+o(1))\mathbb{E}_{\pi,\tilde{\pi}}\left[\exp\left(\frac{\rho^2}{1-\rho^{2}}\left[\binom{n_1}{m}+n_2\binom{n_1}{m-2}\right]\right)\right]\\
&=&(1+o(1))\mathbb{E}_{\pi,\tilde{\pi}}\left[\exp\left(\frac{\rho^2}{1-\rho^{2}}\left[\binom{n_1}{m}+n_2\binom{n_1}{m-2}\right]\right)I[0\leq n_1< \sqrt{n}]\right]\\
&&+(1+o(1))\mathbb{E}_{\pi,\tilde{\pi}}\left[\exp\left(\frac{\rho^2}{1-\rho^{2}}\left[\binom{n_1}{m}+n_2\binom{n_1}{m-2}\right]\right)I[\sqrt{n}\leq n_1\leq n]\right]\\
&=&(a)+(b).
\end{eqnarray*}
If $n_1<\sqrt{n}$, then 
\[\rho^2\left[\binom{n_1}{m}+n_2\binom{n_1}{m-2}\right]=O\left( \frac{n^{\frac{m}{2}}\log n }{n^{m-1}}+\frac{n^{1+\frac{m-2}{2}}\log n }{n^{m-1}}\right)=o(1),\ \ \ m\geq3.\]
Hence $(a)=1+o(1)$.

Next, we show $(b)=o(1)$ if $\rho^2<\frac{(1-\epsilon)n\log n}{\binom{n}{m}}$ . Let $Z_t$, ($1\leq t\leq k$) be independent Poisson variables with $Z_t\sim Poi(\frac{1}{t})$. By Lemma \ref{randomperm}, we have
\begin{eqnarray}\nonumber
(b)&\leq&(1+o(1))\mathbb{E}_{Z_1,Z_2}\left[\exp\left(\frac{\rho^2}{1-\rho^{2}}\left[\binom{Z_1}{m}+Z_2\binom{Z_1}{m-2}\right]\right)I[\sqrt{n}\leq Z_1\leq n]\right]e^{\frac{3}{2}}\\ \label{bz}
&=&(1+o(1))e^{\frac{3}{2}}\mathbb{E}_{Z_1}\left[\exp\left(\frac{\rho^2}{1-\rho^{2}}\binom{Z_1}{m}\right)I[\sqrt{n}\leq Z_1\leq n]\mathbb{E}_{Z_2}\left(\exp\left(\frac{\rho^2}{1-\rho^{2}}Z_2\binom{Z_1}{m-2}\right)\Bigg|Z_1\right)\right].
\end{eqnarray}
By the moment generating function of Poisson distribution, we have
\begin{eqnarray*}
\mathbb{E}_{Z_2}\left(\exp\left(\frac{\rho^2}{1-\rho^{2}}Z_2\binom{Z_1}{m-2}\right)\Bigg|Z_1\right)&=&\exp\left[\frac{1}{2}\left(e^{\frac{\rho^2}{1-\rho^{2}}\binom{Z_1}{m-2}}-1\right)\right].
\end{eqnarray*}
On the event $\sqrt{n}\leq Z_1\leq n$, it follows that
\[\frac{\rho^2}{1-\rho^{2}}\binom{Z_1}{m-2}=O\left(\frac{\log n}{n^{m-1}}n^{m-2}\right)=o(1).\]
Hence, by (\ref{bz}), $k!\geq \left(\frac{k}{e}\right)^k$ and $\rho^2<\frac{(1-\epsilon)n\log n}{\binom{n}{m}}$, we have
\begin{eqnarray}\nonumber
(b)&\leq&(1+o(1))e^{\frac{3}{2}}\mathbb{E}_{Z_1}\left[\exp\left(\frac{\rho^2}{1-\rho^{2}}\binom{Z_1}{m}\right)I[\sqrt{n}\leq Z_1\leq n]\right]\\ \nonumber
&=&(1+o(1))e^{\frac{3}{2}-1}\sum_{k=\sqrt{n}}^n\exp\left(\frac{\rho^2}{1-\rho^{2}}\binom{k}{m}\right)\frac{1}{k!}\\ \nonumber
&\leq&(1+o(1))e^{\frac{3}{2}-1}\sum_{k=\sqrt{n}}^n\exp\left(\frac{\rho^2}{1-\rho^{2}}\binom{k}{m}-k\log k-k\right)\\ \label{be}
&\leq&(1+o(1))e^{\frac{3}{2}-1}\sum_{k=\sqrt{n}}^n\exp\left(k\left((1-\epsilon)\frac{\log n}{n^{m-1}}k^{m-1}-\log k\right)-k\right).
\end{eqnarray}
Define $f(k)=(1-\epsilon)\frac{\log n}{n^{m-1}}k^{m-1}-\log k$. The derivative of $f(k)$ is equal to
\[f^{\prime}(k)=(1-\epsilon)(m-1)\frac{\log n}{n^{m-1}}k^{m-2}-\frac{1}{k}.\]
Solving $f^{\prime}(k)=0$ yields $k_0=\frac{n}{\left((1-\epsilon)(m-1)\log n\right)^{\frac{1}{m-1}}}$. Then $f(k)$ is decreasing for $k\leq k_0$ and increasing for $k\geq k_0$. Hence, 
\begin{eqnarray*}
f(k)\leq\max\left\{f(\sqrt{n}),f(n)\right\}&=& \max\left\{(1-\epsilon)\frac{\log n}{n^{m-1}}\sqrt{n}^{m-1}-\log \sqrt{n},(1-\epsilon)\frac{\log n}{n^{m-1}}n^{m-1}-\log n\right\}\\
&=&\max\left\{-\frac{1}{2}\log n(1+o(1)),-\epsilon\log n\right\}.
\end{eqnarray*}

By (\ref{be}), for a positive consant $c$, we have
\begin{eqnarray}\nonumber
(b)&\leq&(1+o(1))e^{\frac{3}{2}-1}e^{\left(\log n-\sqrt{n}-c\sqrt{n}\log n\right)}=o(1).
\end{eqnarray}
Then the proof is complete.
\end{proof}

The bound in Proposition \ref{prho} is not sharp. The conditional second moment method will be used to close the gap. The result is summarized in the following Proposition \ref{prho2}.

\begin{Proposition}\label{prho2}
 If $\frac{n\log n}{\binom{n}{m}}\leq\rho^2<\frac{(1-\epsilon)2n\log n}{\binom{n}{m}}$ for any positive constant $\epsilon$, then $H_0$ and $H_1$ are indistinguishable.
\end{Proposition}

\begin{proof}[Proof of Proposition \ref{prho2}] 
We use the conditional second moment method as in \cite{WXY21} to prove Proposition \ref{prho2}. 

Let $I$ be the set of fixed points of $\sigma$ and $\mathcal{O}_1$ be the set of subsets in $I$ with cardinality $m$. Then for any $\{i_1,\dots,i_m\}\in \mathcal{O}_1$, $\{\pi_{i_1},\dots,\pi_{i_m}\}=\{\tilde{\pi}_{i_1},\dots,\tilde{\pi}_{i_m}\}$. For $S\subset[n]$ and a positive constant $C$, define event $E_S$ as
\begin{eqnarray*}
E_S&=&\Bigg\{\sum_{\{i_1,i_2,\dots i_m\}\subset S}A_{1,i_1i_2\dots i_m}^2\geq\binom{|S|}{m}-t_S,\sum_{\{i_1,i_2,\dots i_m\}\subset S}A_{2,\pi_{i_1}\dots\pi_{i_m}}^2\geq\binom{|S|}{m}-t_S,\\
&&\sum_{\{i_1,i_2,\dots i_m\}\subset S}A_{1,i_1i_2\dots i_m}A_{2,\pi_{i_1}\dots\pi_{i_m}}\leq\rho\binom{|S|}{m}+t_S\Bigg\},
\end{eqnarray*}
where $t_S$ is of order $n^{\frac{m+1}{2}}$.
Let 
\[E=\cap_{S\subset[n],|S|\geq\frac{n}{2^{\frac{1}{m-1}}}}E_S.\]

By Lemma \ref{lemmap}, $\mathbb{P}(E)=1-o(1)$ under $H_1$. Hence, we have
\begin{equation}\label{es2}
\mathbb{E}\left[\left(\frac{Q(A_1,A_2)}{P(A_1,A_2)}\right)^2\right]=\mathbb{E}_{\pi,\tilde{\pi}}\left(\prod_{O\in\mathcal{O}}\mathbb{E}(L_{O}I[A_1,A_2,\pi\in E]I[A_1,A_2,\tilde{\pi}\in E])\right).
\end{equation}

For $n_1\leq \frac{n}{2^{\frac{1}{m-1}}}$, by a similar argument as in the proof of Lemma \ref{enrho}, one has
\begin{eqnarray*}
&&\prod_{O\in\mathcal{O}}\mathbb{E}(L_{O}I[A_1,A_2,\pi\in E]I[A_1,A_2,\tilde{\pi}\in E])\leq \mathbb{E}\left(\prod_{O\in\mathcal{O}}L_{O}\right)=\prod_{O\in\mathcal{O}}\frac{1}{1-\rho^{2|O|}}\\
&=&\prod_{O\in\mathcal{O}_1}\frac{1}{1-\rho^{2|O|}}\prod_{O\notin\mathcal{O}_1}\frac{1}{1-\rho^{2|O|}}\\
&=&\left(\frac{1}{1-\rho^2}\right)^{\binom{n_1}{m}+n_2\binom{n_1}{m-2}}\left(\frac{1}{1-\rho^4}\right)^{\binom{n}{m}}(1+o(1))\\
&\leq&\exp\left(\frac{\rho^2}{1-\rho^2}\left[\binom{n_1}{m}+n_2\binom{n_1}{m-2}\right]\right)(1+o(1)).
\end{eqnarray*}

Suppose $n_1\geq \frac{n}{2^{\frac{1}{m-1}}}$. Since $\rho^2\geq\frac{n\log n}{\binom{n}{m}}$, then $n^{\frac{m+1}{2}}=o(\rho\binom{n_1}{m})$. In this case, on event $E_I$, we get
\[\sum_{\{i_1,i_2,\dots i_m\}\subset S}A_{1,i_1i_2\dots i_m}^2\geq\binom{|S|}{m}(1+o(1)),\sum_{\{i_1,i_2,\dots i_m\}\subset S}A_{2,\pi_{i_1}\dots\pi_{i_m}}^2\geq\binom{|S|}{m}(1+o(1)),\]
\[\sum_{\{i_1,i_2,\dots i_m\}\subset S}A_{1,i_1i_2\dots i_m}A_{2,\pi_{i_1}\dots\pi_{i_m}}\leq\rho\binom{|S|}{m}(1+o(1)).\]
Then it follows that
\begin{eqnarray*}
&&\mathbb{E}\left[\prod_{O\in\mathcal{O}}L_{O}I[A_1,A_2,\pi\in E]I[A_1,A_2,\tilde{\pi}\in E]\right]\\
&\leq& \mathbb{E}\left[\prod_{O\in\mathcal{O}}L_{O}I[A_1,A_2,\pi\in E_I]\right]\\
&=&\mathbb{E}\left[\prod_{\{i_1,i_2,\dots i_m\}\subset I}L_{i_1i_2\dots i_m}I[A_1,A_2,\pi\in E_I]\right]\prod_{O\notin\mathcal{O}_1}\frac{1}{1-\rho^{2|O|}}.
\end{eqnarray*}
Further, on event $E_I$, the following inequalities hold.
\begin{eqnarray*}
&&\mathbb{E}\left[\prod_{\{i_1,i_2,\dots i_m\}\subset I}L_{i_1i_2\dots i_m}I[A_1,A_2,\pi\in E_I]\right]\\
&\leq&\frac{1}{(1-\rho^2)^{\binom{n_1}{m}}}\mathbb{E}\Bigg[\exp\left\{-\frac{\rho^2}{(1-\rho^2)}\sum_{\{i_1,i_2,\dots i_m\}\subset I}(A_{1,i_1i_2\dots i_m}^2+A_{2,\pi_{i_1}\pi_{i_2}\dots \pi_{i_m}}^2)\right\}\\ \label{likelihood}
&&\times\exp\left\{\frac{2\rho}{1-\rho^2}\sum_{\{i_1,i_2,\dots i_m\}\subset I}A_{1,i_1i_2\dots i_m}A_{2,\pi_{i_1}\pi_{i_2}\dots \pi_{i_m}}\right\}I[A_1,A_2,\pi\in E_I]\Bigg]\\
&\leq&\frac{1}{(1-\rho^2)^{\binom{n_1}{m}}}\exp\left\{-\frac{(1+o(1))2\rho^2}{(1-\rho^2)}\binom{n_1}{m}\right\}\\
&&\times\mathbb{E}\Bigg[\exp\left\{\frac{2\rho}{1-\rho^2}\sum_{\{i_1,i_2,\dots i_m\}\subset I}A_{1,i_1i_2\dots i_m}A_{2,\pi_{i_1}\pi_{i_2}\dots \pi_{i_m}}\right\}I\left[\sum_{\{i_1,i_2,\dots i_m\}\subset I}A_{1,i_1i_2\dots i_m}A_{2,\pi_{i_1}\pi_{i_2}\dots \pi_{i_m}}\leq \rho\binom{n_1}{m}\right]\Bigg]\\
&\leq&\frac{1}{(1-\rho^2)^{\binom{n_1}{m}}}\exp\left\{-\frac{(1+o(1))2\rho^2}{(1-\rho^2)}\binom{n_1}{m}\right\}\\
&&\times\mathbb{E}\Bigg[\exp\left\{\frac{2\rho}{1-\rho^2}\left(\frac{1-\rho^2}{2}\sum_{\{i_1,i_2,\dots i_m\}\subset I}A_{1,i_1i_2\dots i_m}A_{2,\pi_{i_1}\pi_{i_2}\dots \pi_{i_m}}+(1-\frac{1-\rho^2}{2})\rho\binom{n_1}{m}\right)\right\} \Bigg]\\
&=&\frac{1}{(1-\rho^2)^{\binom{n_1}{m}}}\exp\left\{-\frac{(1+o(1))2\rho^2}{(1-\rho^2)}\binom{n_1}{m}\right\}\exp\left\{\frac{\rho^2(1+\rho^2)}{1-\rho^2}\binom{n_1}{m}\right\}\exp\left\{-\frac{1}{2}\binom{n_1}{m}\log(1-\rho^2)\right\} \\
&=&\exp\left\{\frac{(1+o(1))\rho^2}{2}\binom{n_1}{m}\right\}.
\end{eqnarray*}
In the second last equality we used the fact that $\mathbb{E}[e^{\lambda XY}]=\frac{1}{1-\lambda^2}$ for independent standard normal random variables $X,Y$ and $|\lambda|<1$.

Then we can bound the second moment of the likelihood ratio under $H_0$ as
\begin{eqnarray*} 
\mathbb{E}\left[\left(\frac{Q(A_1,A_2)}{P(A_1,A_2)}\right)^2\right]&=&(1+o(1))\mathbb{E}\left[\exp\left(\frac{\rho^2}{1-\rho^2}\left[\binom{n_1}{m}+n_2\binom{n_1}{m-2}\right]\right)I[n_1\leq\frac{n}{2^{\frac{1}{m-1}}}]\right]\\
&&+(1+o(1))\mathbb{E}\left[\exp\left\{\frac{\rho^2}{2}\binom{n_1}{m}+\frac{\rho^2}{1-\rho^2}n_2\binom{n_1}{m-2}\right\}I[n_1\geq\frac{n}{2^{\frac{1}{m-1}}}]\right]\\
&=&(c)+(d).
\end{eqnarray*}

By the proof of Lemma \ref{enrho} and $\rho^2<\frac{(1-\epsilon)2n\log n}{\binom{n}{m}}$, we have
\begin{eqnarray*} 
(c)&=&(1+o(1))
\mathbb{E}\left[\exp\left(\frac{\rho^2}{1-\rho^2}\binom{n_1}{m}\right)I[n_1\leq\sqrt{n}]\right]+(1+o(1))
\mathbb{E}\left[\exp\left(\frac{\rho^2}{1-\rho^2}\binom{n_1}{m}\right)I[\sqrt{n}<n_1\leq\frac{n}{2^{\frac{1}{m-1}}}]\right]\\
&\leq&1+o(1)+e^{\frac{3}{2}-1}\sum_{k=\sqrt{n}}^{\frac{n}{2^{\frac{1}{m-1}}}}\exp\left(k\left((1-\epsilon)\frac{2\log n}{n^{m-1}}k^{m-1}-\log k\right)-k\right).
\end{eqnarray*}
Let $f(k)=(1-\epsilon)\frac{2\log n}{n^{m-1}}k^{m-1}-\log k$. Similar to the proof of Lemma \ref{enrho}, it is easy to verify
\[f(k)\leq\max\left\{f(\sqrt{n}),f\left(\frac{n}{2^{\frac{1}{m-1}}}\right)\right\}=\max\left\{-\frac {1}{2}\log n(1+o(1)),-\epsilon\log n (1+o(1))\right\}.\]
Hence, $(c)=1+o(1)$.

For $(d)$, by Lemma \ref{randomperm}, one has
\begin{eqnarray*} 
(d)&=&(1+o(1))
\mathbb{E}\left[\exp\left(\frac{\rho^2}{2}\binom{n_1}{m}\right)I[n_1\geq\frac{n}{2^{\frac{1}{m-1}}}]\right]\\
&\leq&e^{\frac{3}{2}-1}\sum_{k=\frac{n}{2^{\frac{1}{m-1}}}}^{n}\exp\left(k\left((1-\epsilon)\frac{\log n}{n^{m-1}}k^{m-1}-\log k\right)-k\right).
\end{eqnarray*}
Let $f(k)=(1-\epsilon)\frac{\log n}{n^{m-1}}k^{m-1}-\log k$. Similar to the proof of Lemma \ref{enrho}, it is easy to verify
\[f(k)\leq\max\left\{f(n),f\left(\frac{n}{2^{\frac{1}{m-1}}}\right)\right\}=\max\left\{-\frac {1}{2}\log n(1+o(1)),-\frac{\epsilon+1}{2}\log n\right\}.\]
Hence, $(d)=o(1)$. Then it follows that
\[\mathbb{E}\left[\left(\frac{Q(A_1,A_2)}{P(A_1,A_2)}\right)^2\right]\leq 1+o(1).\]
The proof is complete.
\end{proof}

\begin{Lemma}\label{lemmap}
Under $H_1$,
$\mathbb{P}(E)=1-o(1)$. 
\end{Lemma}
\begin{proof}[Proof of Lemma \ref{lemmap}]
For integer $k$ with $ \frac{n}{2^{\frac{1}{m-1}}}\leq k\leq n$, let $\delta_k=2^{-k}\binom{n}{k}^{-1}$, $S$ be a subset with $|S|=k$ and $t_S=C\left(\sqrt{\binom{k}{m}\log\frac{1}{\delta_k}}+\log\frac{1}{\delta_k}\right)=Cn^{\frac{m+1}{2}}(1+o(1))$. By Hanson-Wright inequality in Lemma \ref{hansonwright}, we have $\mathbb{P}(E_S^c)\leq 6\delta_k$. Hence,
\[\mathbb{P}(E^c)\leq 6\sum_{k=\frac{n}{2^{\frac{1}{m-1}}}}^n\binom{n}{k}\delta_k\leq 6n2^{-\frac{n}{2^{\frac{1}{m-1}}}}=o(1).\]
Then the proof is complete.
\end{proof}

\section{Erd\"{o}s-R\'{e}nyi Hypergraph}

In this section, we study the hypergraph correlation test under the Erd\"{o}s-R\'{e}nyi model. In this case, the hypothesis (\ref{hypothesis}) is reformulated as follows. 
\begin{eqnarray}\nonumber
&&H_0: A_{1,i_1i_2\dots i_m},
A_{2,i_1i_2\dots i_m}\overset{i.i.d.}{\sim}Bern(ps),\\ \nonumber
&&
H_1: A_{1,i_1i_2\dots i_m}\overset{i.i.d.}{\sim}Bern(ps), \ A_{2,\pi_{i_1}\pi_{i_2}\dots \pi_{i_m}}\overset{i.i.d.}{\sim}Bern\left(sA_{1,i_1i_2\dots i_m}+(1-A_{1,i_1i_2\dots i_m})\frac{ps(1-s)}{1-ps}\right),\\ \label{erhypo}
&&\hskip 7cm\ conditional\ on\ \pi\sim Unif(P_n).
 \label{erhyp}
\end{eqnarray}
It is easy to verify the correlation between $A_{1,i_1i_2\dots i_m}$ and $A_{2,\pi_{i_1}\pi_{i_2}\dots \pi_{i_m}}$ under $H_1$ is 
\[\rho=\frac{s(1-p)}{1-ps}.\]
When $p=o(1)$, $\rho=s(1+o(1))$. In this case, $s$ measures the scale of correlation. For $m=2$, the  correlated Erd\"{o}s-R\'{e}nyi graph model is proposed in \cite{PG11} and widely studied in graph matching problem (\cite{BCLSS19,MX19,DMWX21,WXY21,CWXY21}).

The following theorem provides a sharp testing threshold for hypothesis (\ref{erhyp}) when the Erd\"{o}s-R\'{e}nyi hypergraphs are dense.
\begin{Theorem}[Erd\"{o}s-R\'{e}nyi model]\label{thm:2}
Let $m\geq2$ be a fixed integer. Then $H_0$ and $H_1$ in (\ref{erhypo}) are distinguishable if
\[s^2\geq\frac{n\log n}{\binom{n}{m}\left(\log\frac{1}{p}-1+p\right)p}.\]
Suppose $p$ is bounded away from one and $\log\frac{1}{p}=o(\log n)$. Then $H_0$ and $H_1$ in (\ref{erhypo}) are indistinguishable if
\begin{equation}\label{sc5}
s^2<\frac{(1-\epsilon)n\log n}{\binom{n}{m}\left(\log\frac{1}{p}-1+p\right)p},
\end{equation}
for any constant $\epsilon>0$.
\end{Theorem}

 For Erd\"{o}s-R\'{e}nyi model, the sharp testing boundary is $\frac{n\log n}{\binom{n}{m}\left(\log\frac{1}{p}-1+p\right)p}$, which decreases as $m$ gets larger. This shows that testing correlated Erd\"{o}s-R\'{e}nyi hypergraph ($m\geq3$) is easier than testing correlated Erd\"{o}s-R\'{e}nyi graphs (see result for $m=2$ in \cite{WXY21}).


\begin{proof}[Proof of Theorem \ref{thm:2}] 
{\bf (Positive result).} Similar to the Gaussian Wigner model, we shall use the maximum likelihood method to construct a powerful test statistic. The likelihood ratio given $\pi$ is equal to 
\begin{eqnarray*}
\frac{Q(A_1,A_2|\pi)}{P(A_1,A_2)}&=&\prod_{1\leq i_1<\dots<i_m\leq n}(1-s)^{A_{1,i_1i_2\dots i_m}}\left(\frac{1-2ps+ps^2}{1-ps}\right)^{1-A_{1,i_1i_2\dots i_m}}\\
&&\times \prod_{1\leq i_1<\dots<i_m\leq n}\frac{1}{1-ps}\left(\frac{(1-ps)(1-s)}{1-2ps+ps^2}\right)^{A_{2,\pi_{i_1}\pi_{i_2}\dots \pi_{i_m}}}\\
&&\times \left(\frac{1-2ps+ps^2}{p(1-s)^2}\right)^{\sum_{1\leq i_1<\dots<i_m\leq n}A_{1,i_1i_2\dots i_m}A_{2,\pi_{i_1}\pi_{i_2}\dots \pi_{i_m}}}.
\end{eqnarray*}
Let $T_n=\max_{\pi}T(\pi)$ with $T(\pi)=\sum_{1\leq i_1<\dots<i_m\leq n}A_{1,i_1i_2\dots i_m}A_{2,\pi_{i_1}\pi_{i_2}\dots \pi_{i_m}}$.

The correlation coefficient $\rho$ for Erdos-Renyi model is given by
\[\rho=\frac{s(1-p)}{1-ps}.\]
Larger $s$ implies larger correlation $\rho$. Hence, it is easier to test the correlation. Then we can assume
\begin{equation}\label{ssqure}
s^2=\frac{n\log n}{\binom{n}{m}\left(\log\frac{1}{p}-1+p\right)p},
\end{equation}
which implies $p\gg\frac{1}{n^{m-1}}$ and $\binom{n}{m}ps^2\gg n$. Let $t_n=\binom{n}{m}ps^2(1-\tau_n)$ with $ \left(\binom{n}{m}ps^2\right)^{-0.5}\ll\tau_n<1$.

Under $H_1$, we show $\mathbb{P}(T_n\geq t_n)=1+o(1)$. Note that the product $A_{1,i_1i_2\dots i_m}A_{2,\pi_{i_1}\pi_{i_2}\dots \pi_{i_m}}$ are independent and follow Bernoulli($ps^2$). Hence $T(\pi)\sim Binomial(\binom{n}{m},ps^2)$. By Chenorff bound in Lemma \ref{chernoff}, it is easy to get 
\[\mathbb{P}(T_n\leq t_n)\leq \mathbb{P}(T(\pi)\leq t_n)\leq e^{-\frac{\tau_n^2}{2}\binom{n}{m}ps^2}=o(1).\]

Next, we show under $H_0$, $\mathbb{P}(T_n\geq t_n)=o(1)$. In this case, $A_{1,i_1i_2\dots i_m}A_{2,\pi_{i_1}\pi_{i_2}\dots \pi_{i_m}}$ are independent and follow Bernoulli($p^2s^2$). Hence $T(\pi)\sim Binomial(\binom{n}{m},p^2s^2)$. By the multiplicative Chernoff bound in Lemma \ref{chernoff}, we have
\begin{eqnarray*}
\mathbb{P}(T_n\geq t_n)&\leq&n!\mathbb{P}(T(\pi)\geq t_n)\\
&\leq&n!\exp\left(\binom{n}{m}p^2s^2\left[\frac{1-\tau_n}{p}\log\frac{1-\tau_n}{p}+1-\frac{1-\tau_n}{p}\right]\right)\\
&=&n!\exp\left(\binom{n}{m}ps^2(1-\tau_n)\log\frac{1-\tau_n}{ep}-\binom{n}{m}p^2s^2\right)\\
&\leq&n!\exp\left[-\binom{n}{m}ps^2\left(\log\frac{1}{p}-1+p\right)+\tau_n\binom{n}{m}ps^2\log\frac{1}{p}\right]\\
&\leq&e\exp\left[-n+\tau_n\binom{n}{m}ps^2\log\frac{1}{p}+0.5\log n\right].
\end{eqnarray*}
If $p$ is bounded away from one, then $\binom{n}{m}ps^2=O(n\log n)$. Taking $\tau_n=\left(\binom{n}{m}ps^2\right)^{-0.5}\log n$ and noting that $\log \frac{1}{p}=o(\log n)$ yields $\mathbb{P}(T_n\geq t_n)=o(1)$. 

Suppose $p=1+o(1)$. Take $\tau_n=\left(\binom{n}{m}ps^2\right)^{-\epsilon}$ with $\frac{m-1}{2m-1}<\epsilon<0.5$. Note that for some positive constant $c>0$, by (\ref{ssqure}) it follows that
\[\log\frac{1}{p}-1+p=\frac{(1-p)^2}{p}(1-\frac{1}{2p})(1+o(1))\geq c\frac{\log n}{n^{m-1}}.\]
Besides, $\log\frac{1}{p}<\frac{1-p}{p}$. Hence
\[\tau_n\binom{n}{m}ps^2\log\frac{1}{p}=O\left(\frac{(n\log n)^{1-\epsilon}}{(1-p)^{2(1-\epsilon)}}(1-p)\right)=O\left(\log n)^{\epsilon} n^{1+(m-1)-\epsilon(2m-1)}\right)=o(n).\]
Then $\mathbb{P}(T_n\geq t_n)=o(1)$. The proof is complete.
\end{proof}

\begin{proof}[Proof of Theorem \ref{thm:2}] 
{\bf (Negative result).} To prove the negative result, we prove the second moment of the likelihood ratio under $H_0$ converges to one, that is,
\[\mathbb{E}\left[\left(\frac{Q(A_1,A_2)}{P(A_1,A_2)}\right)^2\right]\leq 1+o(1),\]
under $H_0$.

Assume $p\in(0,1-\epsilon_0)$ for a constant $\epsilon_0\in(0,1)$ and $\log\frac{1}{p}=o(\log n)$. Suppose
\begin{equation}\label{sc2}
s^2=\frac{(1-\epsilon)n\log n}{\binom{n}{m}\left(\log\frac{1}{p}-1+p\right)p},
\end{equation}
for any constant $\epsilon\in(0,1)$. In this case,
\begin{equation}\label{nr1}
s^2=o(1),\ \ s^2\gg\frac{n}{\binom{n}{m}}=n^{-(m-1)},\ \ \ \binom{n}{m}ps^2\gg n.
\end{equation}

Let $w(x)$ be the solution of equation $w(x)e^{w(x)}=x$ for $x\geq-\frac{1}{e}$.
Define $\zeta(k)$ as
\[\zeta(k)=\binom{k}{m}ps^2\exp\left(1+w\left(\frac{k\log\frac{2en}{k}}{eps^2\binom{k}{m}}-\frac{1}{e}\right)\right),\ \ k\geq m.\]
Let 
\[\alpha_p=\left(\log\frac{1}{p}-1+p\right)p.\]
Clearly, $\alpha_p\geq cn^{-o(1)}$ for some constant $c>0$ and $n\alpha_p^{\frac{1}{m-1}}\geq cn^{1-o(1)}$. Define event $E$ as
\[E=\cap_{n\alpha_p^{\frac{1}{m-1}}\leq|S|\leq n,S\subset[n]}E_S,\]
where $E_S$ is given by
\begin{eqnarray*}
E_S&=&\Bigg\{\sum_{\{i_1,i_2,\dots i_m\}\subset S}A_{1,i_1i_2\dots i_m}\geq\binom{|S|}{m}ps-\sqrt{2\binom{|S|}{m}ps|S|\log\frac{2en}{|S|}},\\
&&\sum_{\{i_1,i_2,\dots i_m\}\subset S}A_{2,\pi_{i_1}\dots\pi_{i_m}}\geq\binom{|S|}{m}ps-\sqrt{2\binom{|S|}{m}ps|S|\log\frac{2en}{|S|}},\\
&&\sum_{\{i_1,i_2,\dots i_m\}\subset S}A_{1,i_1i_2\dots i_m}A_{2,\pi_{i_1}\dots\pi_{i_m}}\leq\zeta(|S|)
\Bigg\}
\end{eqnarray*}

\begin{Lemma}\label{lnr2}
Under $H_1$, $\mathbb{P}(E)=1-e^{-\Omega(n\alpha_p^{\frac{1}{m-1}})}$.
\end{Lemma}
\begin{proof}[Proof of Lemma \ref{lnr2}] For $S\subset[n]$ with $|S|=k$, let $\delta_k=\left(\frac{k}{2en}\right)^k$ and $t_n=\sqrt{2\binom{|S|}{m}ps\log\frac{1}{\delta_k}}$ and
\[v_n=\binom{k}{m}ps^2\exp\left(1+w\left(\frac{\log\frac{1}{\delta_k}}{eps^2\binom{k}{m}}-\frac{1}{e}\right)\right).\]
By the multiplicative Chernoff bound in Lemma \ref{chernoff}, we have
\[\mathbb{P}\left(\sum_{\{i_1,i_2,\dots i_m\}\subset S}A_{1,i_1i_2\dots i_m}\leq\binom{|S|}{m}ps-t_n\right)\leq\exp\left(-\log\frac{1}{\delta_k}\right)=\left(\frac{k}{2en}\right)^k,\]
\[\mathbb{P}\left(\sum_{\{i_1,i_2,\dots i_m\}\subset S}A_{2,\pi_{i_1}\dots\pi_{i_m}}\leq\binom{|S|}{m}ps-t_n\right)\leq\exp\left(-\log\frac{1}{\delta_k}\right)=\left(\frac{k}{2en}\right)^k,\]
\[\mathbb{P}\left(\sum_{\{i_1,i_2,\dots i_m\}\subset S}A_{1,i_1i_2\dots i_m}A_{2,\pi_{i_1}\dots\pi_{i_m}}\geq v_n\right)\leq\exp\left(-\log\frac{1}{\delta_k}\right)=\left(\frac{k}{2en}\right)^k.\]
Hence,
\begin{eqnarray*}
\mathbb{P}(E^c)\leq \sum_{k=n\alpha_p^{\frac{1}{m-1}}}^n\binom{n}{k}3\delta_k\leq 3\sum_{k=n\alpha_p^{\frac{1}{m-1}}}^n\frac{1}{2^k}=e^{-\Omega(n\alpha_p^{\frac{1}{m-1}})}.
\end{eqnarray*}
Then the proof is complete.
\end{proof}

Let
\[L_{i_1i_2\dots i_m}=L_1(A_{1,i_1i_2\dots i_m},A_{2,\pi_{i_1}\pi_{i_2}\dots \pi_{i_m}})L_1(A_{1,i_1i_2\dots i_m},A_{2,\tilde{\pi}_{i_1}\tilde{\pi}_{i_2}\dots \tilde{\pi}_{i_m}}),\]
with
\[L_1(x,y)=\frac{1-\eta}{1-ps}\left(\frac{1-s}{1-\eta}\right)^{x+y}\left(\frac{s(1-\eta)}{\eta(1-s)}\right)^{xy},\ \ \eta=\frac{ps(1-s)}{1-ps}.\]
Then by Lemma \ref{lnr2} we have
\begin{equation}\label{ns3}
\mathbb{E}\left[\left(\frac{Q(A_1,A_2)}{P(A_1,A_2)}\right)^2\right]=(1+o(1))\mathbb{E}_{\pi,\tilde{\pi}}\left(\prod_{O\in\mathcal{O}}\mathbb{E}(L_{O}I[A_1,A_2,\pi\in E]I[A_1,A_2,\tilde{\pi}\in E])\right).
\end{equation}

If $n_1\leq n\alpha_p^{\frac{1}{m-1}}$, then
\[\prod_{O\in\mathcal{O}}\mathbb{E}(L_{O}I[A_1,A_2,\pi\in E]I[A_1,A_2,\tilde{\pi}\in E])\leq\prod_{O\in\mathcal{O}}\mathbb{E}(L_{O})=\prod_{O\in\mathcal{O}}(1+\rho^{2|O|}).\]

If $n_1>n\alpha_p^{\frac{1}{m-1}}$, then 
\begin{eqnarray*}
&&\prod_{O\in\mathcal{O}}\mathbb{E}(L_{O}I[A_1,A_2,\pi\in E]I[A_1,A_2,\tilde{\pi}\in E])
\leq\prod_{O\in\mathcal{O}}\mathbb{E}(L_{O}I[E_I]))\\
&=&\prod_{O\notin\mathcal{O}_1}\mathbb{E}(L_{O})\prod_{O\in\mathcal{O}_1}\mathbb{E}(L_{O}I[E_I]))\\
&=&\prod_{O\notin\mathcal{O}_1}(1+\rho^{2|O|})\prod_{\{i_1,\dots,i_m\}\subset I}   \mathbb{E}(L_{i_1i_2\dots i_m}I[E_I])
\end{eqnarray*}

Note that for $n_1>n\alpha_p^{\frac{1}{m-1}}\geq cn^{1-o(1)}$,
\[
\frac{t_n^2}{\binom{n_1}{m}^2p^2s^2}=O\left(\frac{\log\frac{n}{n_1}}{n_1^{m-1}ps}\right)=o\left(\frac{\log n}{n^{\frac{m-1}{2}-o(1)}}\right)=o(1).
\]
Hence, on $E_I$ with $n_1>n\alpha_p^{\frac{1}{m-1}}$, one has \[\sum_{\{i_1,i_2,\dots i_m\}\subset I}A_{1,i_1i_2\dots i_m}\geq\binom{n_1}{m}ps(1+o(1)),\ \ \ \sum_{\{i_1,i_2,\dots i_m\}\subset I}A_{2,\pi_{i_1}\dots\pi_{i_m}}\geq\binom{n_1}{m}ps(1+o(1)),\]
\[\sum_{\{i_1,i_2,\dots i_m\}\subset I}A_{1,i_1i_2\dots i_m}A_{2,\pi_{i_1}\dots\pi_{i_m}}\leq\zeta(n_1).\]
Next, we consider the order of $\zeta(n_1)$. Define
\[\gamma=\frac{n_1\log\frac{2en}{n_1}}{\binom{n_1}{m}ps^2}.\]

\begin{Lemma}\label{lnr3}
$[I]$. If $\gamma=o(1)$, then $\zeta(n_1)=\binom{n_1}{m}ps^2(1+o(1))$. \\
$[II]$. If $\gamma=\Theta(1)$, then $\zeta(n_1)=\Theta\left(\binom{n_1}{m}ps^2\right)$. For $n_1>n\alpha_p^{\frac{1}{m-1}}$, $\zeta(n_1)=o\left(\binom{n_1}{m}s^2\right)$. \\
$[III]$. If $\gamma\gg 1$, then $\zeta(n_1)=(e+o(1))\binom{n_1}{m}ps^2\frac{\gamma}{\log \gamma}$. For $n_1>n\alpha_p^{\frac{1}{m-1}}$, $\zeta(n_1)=o\left(\binom{n_1}{m}s^2\right)$.
\end{Lemma}
\begin{proof}[Proof of Lemma \ref{lnr3}]$[I]$ follows from the fact that $w(\frac{\gamma-1}{e})=-1+\sqrt{2\gamma}+O(\gamma)$ if $\gamma=o(1)$.

For $[II]$, if $\gamma=\Theta(1)$, it is obvious that  $\zeta(n_1)=\Theta\left(\binom{n_1}{m}ps^2\right)$. Suppose $n_1>n\alpha_p^{\frac{1}{m-1}}$ and $p\geq c$ for some constant $c>0$. Then $s^2=\Theta\left(\frac{\log n}{n^{m-1}}\right)$, $n_1=\Theta(n)$ and
\[\gamma=O\left(\frac{\log\frac{n}{n_1}}{n^{m-1}s^2}\right)=\frac{O(1)}{\log n}=o(1),\]
which contradicts $\gamma=\Theta(1)$. Hence, $p=o(1)$.

For $[III]$, note that $w(x)=\log x-\log\log x+o(1)$ if $x\gg1$. Then $\zeta(n_1)=(e+o(1))\binom{n_1}{m}ps^2\frac{\gamma}{\log \gamma}$. If $n_1>n\alpha_p^{\frac{1}{m-1}}$, then 
\[p\gamma=O\left(\frac{\log\frac{n}{n_1}}{\alpha_pn^{m-1}s^2}\right)=O\left(\frac{\log\frac{1}{p}}{\log n}\right)=o(1).\]
The proof is complete.
\end{proof}

By Lemma \ref{lnr3}, we get $\zeta(n_1)=\binom{n_1}{m}ps^2+o\left(\binom{n_1}{m}s^2\right)$ for $n_1>n\alpha_p^{\frac{1}{m-1}}$. Then

\begin{eqnarray*}
&&\prod_{\{i_1,\dots,i_m\}\subset I}   \mathbb{E}(L_{i_1i_2\dots i_m}I[E_I])\leq\left(\frac{1-\eta}{1-ps}\right)^{2\binom{n_1}{m}}\left(\frac{1-s}{1-\eta}\right)^{4(1+o(1))\binom{n_1}{m}ps}\\
&&\times\mathbb{E}\left[\left(\frac{s(1-\eta)}{\eta(1-s)}\right)^{2\sum_{\{i_1,i_2,\dots i_m\}\subset I}A_{1,i_1i_2\dots i_m}A_{2,\pi_{i_1}\dots\pi_{i_m}}}I[\sum_{\{i_1,i_2,\dots i_m\}\subset I}A_{1,i_1i_2\dots i_m}A_{2,\pi_{i_1}\dots\pi_{i_m}}\leq \zeta(n_1)]\right]\\
&=&\exp\left(-2(1+o(1))\binom{n_1}{m}ps^2(1-p)\right)\\
&&\times \mathbb{E}\left[\lambda^{\sum_{\{i_1,i_2,\dots i_m\}\subset I}A_{1,i_1i_2\dots i_m}A_{2,\pi_{i_1}\dots\pi_{i_m}}}I[\sum_{\{i_1,i_2,\dots i_m\}\subset I}A_{1,i_1i_2\dots i_m}A_{2,\pi_{i_1}\dots\pi_{i_m}}\leq \zeta(n_1)]\right],
\end{eqnarray*}
where $\lambda=\left(\frac{s(1-\eta)}{\eta (1-s)}\right)^{2}=(1+o(1))\frac{1}{p^2}$.

Note that for any $t\in[0,1]$,
\begin{eqnarray*}
&&\mathbb{E}\left[\lambda^{\sum_{\{i_1,i_2,\dots i_m\}\subset I}A_{1,i_1i_2\dots i_m}A_{2,\pi_{i_1}\dots\pi_{i_m}}}I[\sum_{\{i_1,i_2,\dots i_m\}\subset I}A_{1,i_1i_2\dots i_m}A_{2,\pi_{i_1}\dots\pi_{i_m}}\leq \zeta(n_1)]\right]\\
&\leq&\mathbb{E}\left[\lambda^{t\sum_{\{i_1,i_2,\dots i_m\}\subset I}A_{1,i_1i_2\dots i_m}A_{2,\pi_{i_1}\dots\pi_{i_m}}+(1-t)\zeta(n_1)}\right]\\
&=&\lambda^{\zeta(n_1)}\frac{\left(1+(\lambda^t-1)p^2s^2\right)^{\binom{n_1}{m}}}{\lambda^{t\zeta(n_1)}}.
\end{eqnarray*}
Let $g(y)=\frac{\left(1+(y-1)p^2s^2\right)^{\binom{n_1}{m}}}{y^{\zeta(n_1)}}$. It is easy to verify that $g(y)$ attains minimum value at $y_0=\frac{\zeta(n_1)(1-p^2s^2)}{p^2s^2\left(\binom{n_1}{m}-\zeta(n_1)\right)}\in [1,\lambda]$. Let $h(x)=-x\log x-(1-x)\log(1-x)$. Then
\begin{eqnarray*}
\prod_{\{i_1,\dots,i_m\}\subset I}   \mathbb{E}(L_{i_1i_2\dots i_m}I[E_I])&\leq&
\exp\left(-2(1+o(1))\binom{n_1}{m}ps^2(1-p)\right)\lambda^{\zeta(n_1)}\frac{\left(1+(y_0-1)p^2s^2\right)^{\binom{n_1}{m}}}{y_0^{\zeta(n_1)}}\\
&=&\exp\left(-\binom{n_1}{m}ps^2(2-p)+\zeta(n_1)(\log s^2+o(1))+\binom{n_1}{m}h\left(\frac{\zeta(n_1)}{\binom{n_1}{m}}\right)\right)\\
&=&\exp\left(-\binom{n_1}{m}ps^2(2-p)+\zeta(n_1)\log\frac{e\binom{n_1}{m}s^2}{\zeta(n_1)}+o(\zeta(n_1))\right).
\end{eqnarray*}

Then the second moment of the likelihood ratio is bounded by
\begin{eqnarray*}\label{ns3}
&&\mathbb{E}\left[\left(\frac{Q(A_1,A_2)}{P(A_1,A_2)}\right)^2\right]\leq(1+o(1))\mathbb{E}\left[\prod_{O\in\mathcal{O}}(1+\rho^{2|O|})I[n_1\leq n\alpha_p^{\frac{1}{m-1}}]\right]\\
&&+\mathbb{E}\left[\prod_{O\notin\mathcal{O}_1}(1+\rho^{2|O|})\exp\left(-\binom{n_1}{m}ps^2(2-p)+\zeta(n_1)\log\frac{e\binom{n_1}{m}s^2}{\zeta(n_1)}+o(\zeta(n_1))\right)I[n_1\geq n\alpha_p^{\frac{1}{m-1}}]\right].
\end{eqnarray*}

Recall that $\rho=(1+o(1))s$ and $\log\frac{1}{p}=o(\log n)$. Then
\[\rho^4n^m=\frac{\log ^2n}{n^{m-2+o(1)}\log^2\frac{1}{p}}=o(1).\]
By the proof of Lemma \ref{enrho}, we have
\[\prod_{O\notin\mathcal{O}_1}(1+\rho^{2|O|})=(1+\rho^2)^{n_2\binom{n_1}{m-2}(1+o(1))}(1+\rho^{4})^{n^{m}}=(1+o(1))(1+\rho^2)^{n_2\binom{n_1}{m-2}}\leq (1+o(1))e^{\rho^2n_2\binom{n_1}{m-2}}.\]
Besides,
\[\prod_{O\in\mathcal{O}_1}(1+\rho^{2|O|})=(1+\rho^2)^{\binom{n_1}{m}}\leq e^{\rho^2\binom{n_1}{m}}.\]

Then we get
\begin{eqnarray*}\label{ns3}
&&\mathbb{E}\left[\left(\frac{Q(A_1,A_2)}{P(A_1,A_2)}\right)^2\right]\leq(1+o(1))\mathbb{E}\left[e^{\rho^2\binom{n_1}{m}+\rho^2n_2\binom{n_1}{m-2}}I[n_1\leq n\alpha_p^{\frac{1}{m-1}}]\right]\\
&&+(1+o(1))\mathbb{E}\left[\exp\left(\rho^2n_2\binom{n_1}{m-2}-\binom{n_1}{m}ps^2(2-p)+\zeta(n_1)\log\frac{e\binom{n_1}{m}s^2}{\zeta(n_1)}+o(\zeta(n_1))\right)I[n_1\geq n\alpha_p^{\frac{1}{m-1}}]\right]\\
&=&(e)+(f).
\end{eqnarray*}

Next, we are going to show $(e)=1+o(1)$ and $(f)=o(1)$.

We show $(e)=1+o(1)$ first.  Similar to the proof of Lemma \ref{enrho}, we have
\begin{eqnarray*}
&&\mathbb{E}\left[e^{\rho^2\binom{n_1}{m}+\rho^2n_2\binom{n_1}{m-2}}I[n_1\leq n\alpha_p^{\frac{1}{m-1}}]\right]\\
&=&\mathbb{E}\left[e^{\rho^2\binom{n_1}{m}+\rho^2n_2\binom{n_1}{m-2}}I[n_1\leq \sqrt{n}]\right]+\mathbb{E}\left[e^{\rho^2\binom{n_1}{m}+\rho^2n_2\binom{n_1}{m-2}}I[\sqrt{n}<n_1\leq n\alpha_p^{\frac{1}{m-1}}]\right]\\
&\leq&1+o(1)+e^{\frac{3}{2}}\mathbb{E}\left[e^{\rho^2\binom{Z_1}{m}}I[\sqrt{n}<Z_1\leq n\alpha_p^{\frac{1}{m-1}}]\right]\\
&\leq&1+o(1)+e^{-1}\sum_{k=\sqrt{n}}^{n\alpha_p^{\frac{1}{m-1}}}e^{\rho^2\binom{k}{m}-k\log k-k}\\
&=&1+o(1)+e^{-1}\sum_{k=\sqrt{n}}^{n\alpha_p^{\frac{1}{m-1}}}e^{k\left(\frac{(1-\epsilon)\log n}{\alpha_p}\frac{k^{m-1}}{n^{m-1}}-\log k\right)-k}
\end{eqnarray*}
Let $f(k)=\frac{(1-\epsilon)\log n}{\alpha_p}\frac{k^{m-1}}{n^{m-1}}-\log k$. It is easy to see
\[f(k)\leq\max\{f(\sqrt{n}),f(n\alpha_p^{\frac{1}{m-1}})\}=-\min\{0.5,\epsilon\}\log n.\]
Hence $(e)\leq1+o(1)$.

Next, we prove $(f)=o(1)$. To this end, for a large positive constant $C$, define
\[\beta_1=\left(\frac{\log^2\frac{m!\binom{n}{m}ps^2}{n}}{\frac{m!\binom{n}{m}ps^2}{n}}\right)^{\frac{1}{m-1}},\ \ \ \beta_2=\left(\frac{\log\frac{m!\binom{n}{m}ps^2}{n}}{C\frac{m!\binom{n}{m}ps^2}{n}}\right)^{\frac{1}{m-1}}.\]
Then $n_1$ falls in one of the three intervals $[\beta_1n,n]$, $[\beta_2n,\beta_1n]$ and $[n\alpha_p^{\frac{1}{m-1}},\beta_2n]$.

If $n_1\in [\beta_1n,n]$, then $\gamma=o(1)$ and hence $\zeta(n_1)=(1+o(1))\binom{n_1}{m}ps^2$.
In this case,
\begin{eqnarray*}
&&\mathbb{E}\left[\exp\left(\rho^2n_2\binom{n_1}{m-2}-\binom{n_1}{m}ps^2(2-p)+\zeta(n_1)\log\frac{e\binom{n_1}{m}s^2}{\zeta(n_1)}+o(\zeta(n_1))\right)I[n_1\geq n\beta_1]\right]\\
&=&\mathbb{E}\left[\exp\left(\rho^2n_2\binom{n_1}{m-2}-\binom{n_1}{m}ps^2(2-p)+\binom{n_1}{m}ps^2\log\frac{e}{p}\right)I[n_1\geq n\beta_1]\right]\\
&\leq&e^{1.5}\mathbb{E}\left[\exp\left(\binom{Z_1}{m}s^2\alpha_p\right)I[Z_1\geq n\beta_1]\right]\\
&\leq&e^{0.5}\sum_{k=\beta_1n}^ne^{k\left(\frac{(1-\epsilon)\log n}{n^{m-1}}k^{m-1}-\log k\right)-k}=o(1).
\end{eqnarray*}

If $n_1\in[\beta_2n,\beta_1n]$, then $\gamma=\frac{C}{m!(m-1)}(1+o(1))$ and hence $\zeta(n_1)=\Theta\left( \binom{n_1}{m}ps^2\right)$. In this case, for some constant $C_1$,
\begin{eqnarray*}
&&\mathbb{E}\left[\exp\left(\rho^2n_2\binom{n_1}{m-2}-\binom{n_1}{m}ps^2(2-p)+\zeta(n_1)\log\frac{e\binom{n_1}{m}s^2}{\zeta(n_1)}+o(\zeta(n_1))\right)I[\beta_2n\leq n_1\leq n\beta_1]\right]\\
&\leq&\mathbb{E}\left[\exp\left(\rho^2n_2\binom{n_1}{m-2}+C_1\binom{n_1}{m}ps^2\log\frac{e}{p}\right)I[\beta_2n\leq n_1\leq n\beta_1]\right]\\
&\leq&e^{1.5}\mathbb{E}\left[\exp\left(C_1\binom{Z_1}{m}ps^2\log\frac{1}{p}\right)I[\beta_2n\leq Z_1\leq n\beta_1]\right]\\
&\leq&e^{0.5}\sum_{k=\beta_2n}^{\beta_1n}e^{k\left(C_1\frac{(1-\epsilon)\log n}{n^{m-1}}k^{m-1}-\log k\right)-k}=o(1).
\end{eqnarray*}

If $n_1\in[n\alpha_p^{\frac{1}{m-1}},\beta_2n]$, then $\gamma\geq\frac{Cm!}{m-1}(1+o(1))$. For sufficiently large $C$, we have $\zeta(n_1)=O\left(\frac{n_1\log\frac{2en}{n_1}}{\log\gamma}\right)$. Then 
\[
\zeta(n_1)\log\frac{e\binom{n_1}{m}s^2}{\zeta(n_1)}+o(\zeta(n_1))\leq n_1R_n,\]
where
\[R_n=C_2\frac{\log\frac{n}{n_1}}{\log \gamma}\log\frac{n_1^{m-1}s^2\log\gamma}{\log\frac{n}{n_1}},\]
for a constant $C_2$. Suppose $R_n=o(\log n)$, then
\begin{eqnarray*}
&&\mathbb{E}\left[\exp\left(\rho^2n_2\binom{n_1}{m-2}-\binom{n_1}{m}ps^2(2-p)+\zeta(n_1)\log\frac{e\binom{n_1}{m}s^2}{\zeta(n_1)}+o(\zeta(n_1))\right)I[\alpha_p^{\frac{1}{m-1}}n\leq n_1\leq n\beta_2]\right]\\
&\leq&\mathbb{E}\left[\exp\left(\rho^2n_2\binom{n_1}{m-2}-\binom{n_1}{m}ps^2(2-p)+n_1R_n\right)I[\alpha_p^{\frac{1}{m-1}}n\leq n_1\leq n\beta_2]\right]\\
&\leq&e^{1.5}\mathbb{E}\left[\exp\left(-\binom{Z_1}{m}ps^2+Z_1o(\log n)\right)I[\alpha_p^{\frac{1}{m-1}}n\leq n_1\leq n\beta_2]\right]\\
&\leq&e^{0.5}\sum_{k=\alpha_p^{\frac{1}{m-1}}n}^{\beta_2n}e^{-k\left(\frac{(1-\epsilon)\log n}{n^{m-1}(\log\frac{1}{p}-1+p)}k^{m-1}+\log k-o(\log n)\right)-k}=o(1).
\end{eqnarray*}
Here in the last equality we used the fact that $\log\frac{1}{p}-1+p\geq c>0$ for some constant $c$, since $p$ is bounded away from one.

Below we prove $R_n=o(\log n)$. Note that $\log\frac{1}{\alpha_p}=o(\log n)$ and hence
\[
R_n=C_2\frac{\log\frac{n}{n_1}}{\log \gamma}\log\frac{n_1^{m-1}s^2}{\log\frac{n}{n_1}}+C_2\frac{\log\frac{n}{n_1}}{\log \gamma}\log\log\gamma=C_2\frac{\log\frac{n}{n_1}}{\log \gamma}\log\frac{n_1^{m-1}s^2}{\log\frac{n}{n_1}}+o(\log n).
\]
Then it suffices to show 
\begin{equation}\label{rn1}
\frac{\log\frac{n}{n_1}}{\log \gamma}\log\frac{n_1^{m-1}s^2}{\log\frac{n}{n_1}}=o(\log n).
\end{equation}

Let $x=\frac{n}{n_1}\in[\frac{1}{\beta_2},\alpha_p^{-\frac{1}{m-1}}]$. Then $\gamma=\Theta\left(\frac{x^{m-1}\log x}{n^{m-1}ps^2}\right)$. To prove (\ref{rn1}), we only need to prove
\begin{equation}\label{rn2}
\max_{x\in[\frac{1}{\beta_2},\alpha_p^{-\frac{1}{m-1}}]}\frac{\log x}{\log\frac{x^{m-1}\log x}{n^{m-1}ps^2}}\log\frac{n^{m-1}s^2}{x^{m-1}\log n}=o(\log n).
\end{equation}
Let $\delta=\left(\log\frac{\log n}{\log\frac{1}{\alpha_p}}\right)^{-1}$. Then $\delta=o(1)$. To prove (\ref{rn2}), it suffices to show $\psi(x)\leq 0$, with $\psi(x)$ given by
\[\psi(x)=\log (x)\log\frac{n^{m-1}s^2}{x^{m-1}\log n}-\delta \log (n)\log\frac{x^{m-1}\log x}{n^{m-1}ps^2}.\]
Straightforward calculation yields
\[\psi^{\prime}(x)=\frac{\log(n^{m-1}s^2)-2(m-1)\log x-\log\log x-1-(m-1)\delta\log n-\frac{\delta\log n}{\log x}}{x}.\]
It is easy to see that
\[\frac{\log\frac{1}{\alpha_p}}{\delta\log n}=\frac{\log\frac{\log n}{\log\frac{1}{\alpha_p}}}{\frac{\log n}{\log\frac{1}{\alpha_p}}}=o(1).\]
Hence 
\[\log(n^{m-1}s^2)=\log\log n+\log\frac{1}{\alpha_p}=o(\delta\log n).\]
This implies $\psi^{\prime}(x)\leq0$ for $x\in[\frac{1}{\beta_2},\alpha_p^{-\frac{1}{m-1}}]$. Then $\psi(x)\leq\psi(\frac{1}{\beta_2})$. Since
$\frac{1}{\beta_2^{m-1}}\log\frac{1}{\beta_2}=\frac{C}{m-1}(1+o(1))n^{m-1}ps^2$, then

\begin{eqnarray*}
\psi\left(\frac{1}{\beta_2}\right)&=&\log \left(\frac{1}{\beta_2}\right)\log\frac{n^{m-1}s^2}{\frac{1}{\beta_2^{m-1}}\log n}-\delta \log (n)\log\frac{\frac{1}{\beta_2^{m-1}}\log \frac{1}{\beta_2}}{n^{m-1}ps^2}\\
&=&\log \left(\frac{1}{\beta_2}\right)\log\frac{1}{Cp}-\log (C)\delta\log n,
\end{eqnarray*}
which is negative if $p$ is bounded away from zero. Assume $p=o(1)$. Then
\[\frac{1}{\delta}\log \left(\frac{1}{\beta_2}\right)=O\left(\log\frac{\log n}{\log\frac{1}{\alpha_p}}\log\frac{\frac{\log n}{\log\frac{1}{\alpha_p}}}{\log \frac{\log n}{\log\frac{1}{\alpha_p}}}\right)=O\left(\log^2\frac{\log n}{\log\frac{1}{\alpha_p}}\right)=o\left(\frac{\log n}{\log\frac{1}{\alpha_p}}\right),\]
which implies $\psi\left(\frac{1}{\beta_2}\right)\leq0$ for large $n$. Then the proof is complete.
\end{proof}

\section{Additional Lemmas}

In this section, several lemmas are given. Firstly, we present the Hanson-Wright inequality (\cite{WXY21}) below.
\begin{Lemma}[Hanson-Wright]\label{hansonwright}
Let $X,Y\in\mathbb{R}^d$ be standard Gaussion random variables such that $(X_i,Y_i)$, $i=1,2,\dots,d$ are independent and have correlation coefficient $\rho$. Then with probability at $1-2\delta$,
\[|X^TY-d\rho|\leq C\left(d\sqrt{\log\frac{1}{\delta}}+\log\frac{1}{\delta}\right),\]
for a constant $C>0$.
\end{Lemma}

The following lemma presents the Chernoff bound for binomial distribution (\cite{WXY21}).
\begin{Lemma}[Chernoff bound]\label{chernoff}
Let $X\sim Bin(n,p)$ and $\mu=np$. Then for any $\delta>0$,
\[\mathbb{P}(X\geq(1+\delta)\mu)\leq e^{-\mu(1+\delta)\log(1+\delta)-\delta},\]
\[\mathbb{P}(X\leq(1-\delta)\mu)\leq e^{-\frac{\delta}{2}\mu}.\]
Particularly, for $\tau=\mu \exp\left(1+W(\frac{t}{e\mu}-\frac{1}{\mu})\right)$ with $W(x)$ be the solution to the equation $f(x)e^{f(x)}=x$, then
\[\mathbb{P}(X\leq\tau)\leq e^{-t}.\]
\end{Lemma}

The following lemma presents a fact about random permutation (\cite{AT92}).
\begin{Lemma}[]\label{randomperm}
Let $n_t$ be the number of $t$-cycles in a random permutation $\sigma\in P_n$. Let $Z_t\sim Poisson(\frac{1}{t})$ be independent Poisson random variables.. Then
\[\mathbb{E}(g(n_1,n_2,\dots,n_L))\leq e^{1+\frac{1}{2}+\dots+\frac{1}{L}}\mathbb{E}(g(Z_1,Z_2,\dots,Z_L)),\]
for any nonnegative function $g$.
\end{Lemma}

\end{document}